\documentclass[11pt]{article}
\usepackage{amsmath}
\usepackage{graphicx,amssymb,booktabs,multirow}
\usepackage[paperwidth=90mm]{geometry}
\usepackage{anysize}
\usepackage{bm}
\usepackage{natbib,amsmath}
\usepackage{graphicx}
\usepackage{fancyhdr}
\usepackage{rotating}
\marginsize{1in}{1in}{1in}{1in}
\makeatletter
\newcommand*{\indep}{%
	\mathbin{%
		\mathpalette{\@indep}{}%
	}%
}
\newcommand*{\nindep}{%
	\mathbin{
		\mathpalette{\@indep}{\not}
	}%
}
\newcommand*{\@indep}[2]{%
	\sbox0{$#1\perp\m@th$}
	\sbox2{$#1=$}
	\sbox4{$#1\vcenter{}$}
	\rlap{\copy0}
	\dimen@=\dimexpr\ht2-\ht4-.2pt\relax
	\kern\dimen@
	{#2}%
	\kern\dimen@
	\copy0 
}
\makeatother
\usepackage{color}
\usepackage{amsthm}

\makeatother
\newtheorem{theorem}{Theorem}
\newtheorem{lemma}{Lemma}

\newtheorem{assumption}{Assumption}
\newtheorem{corollary}{Corollary}
\usepackage{bbm}
\newcommand{\argmin}{\arg\!\min}

\newcommand{\R}{{\mathbb{R}}}

\usepackage[all]{xy}
\usepackage{bm}
\usepackage{pifont}
\usepackage{stmaryrd}
\usepackage{pgf,tikz}
\usepackage{tikz}
\usetikzlibrary{shapes,arrows,positioning,calc}

\newenvironment{sequation*}{\begin{equation*}\small}{\end{equation*}}

\newenvironment{tequation*}{\begin{equation*}\tiny}{\end{equation*}}

\usepackage{mathrsfs}

\newtheorem{innercustomgeneric}{\customgenericname}
\providecommand{\customgenericname}{}
\newcommand{\newcustomtheorem}[2]{%
	\newenvironment{#1}[1]
	{%
		\renewcommand\customgenericname{#2}%
		\renewcommand\theinnercustomgeneric{##1}%
		\innercustomgeneric
	}
	{\endinnercustomgeneric}
}

\newcustomtheorem{customthm}{Theorem}
\newcustomtheorem{customlemma}{Lemma}
\newcustomtheorem{customassumption}{Assumption}
\tikzstyle{sum} = [draw, circle, minimum size=0pt, inner sep = 1.7pt]

\makeatletter
\newcommand{\mathleft}{\@fleqntrue\@mathmargin0pt}
\newcommand{\mathcenter}{\@fleqnfalse}
\makeatother
\newtheorem{example}{Example}
\newtheorem{definition}{Definition}

\usepackage{algorithm}
\usepackage{algorithmic}
\usepackage{wrapfig}

\allowdisplaybreaks

\usepackage{accents}

\def\T{{ \mathrm{\scriptscriptstyle T} }}

{
	\theoremstyle{plain}

}

\usepackage{subfigure}
\graphicspath{ {./pictures/} }
\usepackage{multirow}

\usepackage{booktabs,caption,fixltx2e}
\usepackage[flushleft]{threeparttable}
\newcommand{\rep}{{\textnormal{rep}}}
\newcommand{\repa}{{\textnormal{repa}}}
\usepackage{enumitem}

\newcommand\Bern{{\textnormal{Bern}}}
\newcommand{\Be}{\textnormal{Be}}
\newcommand{\REPA}{{\texttt{REP-DB}}}
\newcommand{\REP}{{\texttt{REP}}}
\newcommand{\IPW}{{\texttt{IPW}}}
\newcommand{\IPWt}{{\texttt{IPW-true}}}
\newcommand{\IPWu}{{\texttt{IPW-uniform}}}

\newcommand{\bphi}{\bm \phi}
\newcommand{\bPhi}{\bm \Phi}
\newcommand{\bLambda}{\bm \Lambda}
\newcommand{\calN}{\mathcal{N}}
\newcommand{\calR}{\mathcal{R}}

\def\np{\textnormal{np}}

\newtheorem{manualcondinner}{Condition}
\newenvironment{manualcond}[1]{%
	\renewcommand\themanualcondinner{#1}%
	\manualcondinner
}{\endmanualcondinner}

\begin{document}

	\def\spacingset#1{\renewcommand{\baselinestretch}%
		{#1}\small\normalsize} \spacingset{1}

	
	\title{\bf 	Nonparametric inference about mean functionals of nonignorable nonresponse data without identifying the joint distribution}
	
	\author{Wei Li$^1$, Wang Miao$^2$, and Eric Tchetgen Tchetgen$^3$	\vspace{0.5cm}
		\\
		$^1$Center for Applied Statistics and School of Statistics, \\
		Renmin University of China\vspace{0.1cm}
		\\
		$^2$School of Mathematical Sciences, Peking University\\
	$^3$Department of Statistics,
	University of Pennsylvania}
	\date{}
	\maketitle
	
	
	\bigskip
	\begin{abstract}
		We consider  identification and   inference about  mean functionals of observed covariates and an  outcome variable subject to nonignorable missingness.
	By leveraging a  shadow variable, we establish a necessary and sufficient condition for identification of the mean functional even if the full data distribution is not identified. We further  characterize a necessary condition for $\sqrt{n}$-estimability of the mean functional. This condition naturally strengthens the identifying condition, and it
	requires the existence of a function as a solution to a representer equation that connects the shadow variable to the mean functional.
	Solutions to the representer equation may not be unique, which presents   substantial challenges for nonparametric estimation and standard theories for nonparametric sieve estimators are not applicable here.  
	We   construct a consistent estimator for the solution set and then adapt the theory of extremum estimators to find from the estimated set a consistent estimator  for an appropriately chosen solution.
	The estimator is asymptotically normal, locally efficient and attains the semiparametric efficiency bound  under certain regularity conditions.
	We  illustrate the proposed approach via simulations  and a real data application on home pricing. 
	\end{abstract}
	
	\noindent%
	{\it Keywords:}  Identification; Model-free estimation; Nonignorable missingness;  Shadow variable. 
	

	\section{Introduction}

	Nonresponse is  frequently encountered in social science and biomedical studies, 
	due to  such as reluctance to answer sensitive survey questions. 
	Certain characteristics of the missing data mechanism is used to define     a taxonomy to describe the missingness process \citep{rubin1976inference,little2002statistical}.
	It is called missing at random (MAR) if the propensity of missingness conditional on all study variables is unrelated to the missing values. Otherwise, it is called missing not at random (MNAR) or nonignorable.  
	MAR has been commonly  used  for statistical analysis in the presence of missing data;
	however, in many fields of study, suspicion that the  missing data mechanism may be  nonignorable is often warranted   \citep{scharfstein1999adjusting,robins2000sensitivity,rotnitzky1997analysis,rotnitzky1998semiparametric,ibrahim1999missing}.
	For example, nonresponse rates  in surveys about income  tend to be higher for low socio-economic groups \citep{kim2011semiparametric}. In another example, efforts to estimate HIV prevalence in developing countries via household HIV survey and testing such as the well-known Demographic and Health Survey, are likewise subject to nonignorable missing data on participants' HIV status due to highly selective non-participation in the HIV testing component of the survey study \citep{tchetgen2017general}. 
	There  currently exist a variety of methods for the analysis of MAR data, such as likelihood based inference, multiple imputation, inverse probability weighting, and
	doubly  robust methods.
	However, these methods can result in severe  bias and invalid inference in the presence of  nonignorable missing data.

	In this paper, we focus on identification and estimation of mean functionals of observed covariates and 
	an outcome variable subject to nonignorable missingness. Estimating mean functionals is a goal common in many scientific areas, e.g., sampling survey and causal inference, and thus is of significant practical importance.
	However, there are several difficulties for analysis of nonignorable missing data.
	The first challenge is identification, which means that the parameter of interest is uniquely determined from observed data distribution. 
	Identification is straightforward under MAR  as the conditional outcome  distribution in complete-cases equals that in incomplete cases given fully observed covariates,
	whereas it becomes difficult  under MNAR because the selection bias due to missing values is no longer negligible.
	Even if  stringent fully-parametric models are imposed on both the propensity score and the outcome regression, identification may not be achieved; for counterexamples, see \cite{wang2014instrumental,miao2016identifiability}.  
	To resolve the identification difficulty, previous researchers  \citep{robins2000sensitivity,kim2011semiparametric}  have assumed that the selection bias is known or can be estimated from external studies, 
	but this approach should be used rather as a sensitivity analysis, until the validity of the selection bias assumption is assessed.
	Without knowing the selection bias, identification can be achieved by  leveraging   fully observed auxiliary variables that are  available in many empirical studies.
	For instance,  instrumental variables, which are related to the nonresponse propensity but not related to the outcome given covariates, have been used in missing data analysis since \cite{heckman1979sample}. 
	The corresponding semiparametric  theory and inference are recently established by \cite{liu2020identification,sun2018semiparametric,tchetgen2017general,das2003nonparametric}.
	Recently,
	an alternative approach called the shadow variable approach has grown in popularity in sampling
	survey and missing data analysis. In contrast to the instrumental variable, this approach
	entails a shadow variable that is associated with the outcome but independent of the missingness
	process given covariates and the outcome itself.
	Shadow variable 
	is available in many applications \citep{kott2014calibration}. For example, \cite{zahner1992children} and \cite{ibrahim2001using} considered a study of the children's mental health evaluated through their teachers' assessments
	in Connecticut. However, the data for the teachers' assessments are subject to nonignorable missingness. As a proxy of the teacher's assessment, a separate parent report is available for all children in this study. The parent report is likely to be correlated with the teacher's assessment, but is unlikely to be related to the teacher's response rate given the teacher's assessment and fully observed covariates.
	Hence, the parental assessment is regarded as a shadow variable in this study. The shadow variable design is quite general.
	In health and social sciences, an accurately measured outcome is routinely available only for a subset of patients, but one or more surrogates may be fully observed. For instance, \cite{robins1994estimation} considered a cardiovascular disease setting where, due to high cost of laboratory analyses, and to the small amount of stored serum per subject (about $2\%$ of study subjects had stored serum thawed and assayed for antioxidants serum vitamin A and vitamin E), error prone surrogate measurements of the biomarkers derived from self-reported dietary questionnaire were obtained for all subjects. 
	Other important  settings include the semi-supervised set up in comparative effectiveness research where the true outcome is measured only for a small fraction of the data, e.g. diagnosis requiring a costly panel of physicians while surrogates are obtained from databases including ICD-9 codes for certain comorbidities. Instead of assuming MAR conditional on the surrogates, the shadow variable assumption may be more appropriate in presence of informative selection bias.
	By leveraging a shadow variable,   \cite{d2010new,miao2019identification} established  identification results for   nonparametric models under completeness conditions. 
	In related work, \cite{wang2014instrumental,shao2016semiparametric,tang2003analysis,zhao2015semiparametric,morikawa2021semiparametric}
	proposed identification conditions  for a suite of   parametric and  semiparametric models  that require either the propensity score or the  outcome regression, or both to be  parametric. 
	All existing shadow variable approaches need to further impose sufficiently strong conditions to identify the full data distribution, although in practice one may  only be interested in   a  parameter or a mean functional which may be identifiable even if the full data distribution is not.

	The second challenge for MNAR data is the threat of bias due to model misspecification in  estimation, after identification is established.
	Likelihood-based inference 	\citep{greenlees1982imputation,tang2003analysis,tang2014empirical},
	imputation-based methods \citep{kim2011semiparametric},
	inverse probability weighting \citep{wang2014instrumental}
	have been developed for analysis of MNAR data.
	These   estimation methods require correct model specification of either the propensity score or the outcome regression, or both. 
	However,  bias may arise due to specification error of parametric models  as they have limited flexibility, 
	and moreover, model misspecification is more likely  to appear in the presence of missing values.
	\cite{vansteelandt2007estimation,miao2016varieties,miao2019identification,liu2020identification} proposed doubly robust   methods for estimation with MNAR data, which  affords double protection against model misspecification; however, these previous proposals require an odds ratio model characterizing the degree of nonignorable missingness to be either completely known or  correctly specified.

	In this paper, we develop a novel strategy to nonparametrically identify and estimate a generic mean functional. In contrast to previous approaches in the literature, this work has the following distinctive features and makes several contributions to nonignorable missing data literature. 
	First, given a shadow variable, we  directly work on the identification of the mean functional without identifying  the full  data distribution, whereas   previous proposals typically have to first ensure identification of the full data distribution. In particular, we establish a necessary and sufficient  condition for identification of the mean functional, which is weaker than the commonly-used completeness condition for nonparametric identification. For estimation, we propose a representer assumption which is shown to be necessary for $\sqrt{n}$-estimability of the mean functional. The representer assumption
	involves the existence of a function as a
	solution to a representer equation that relates the  mean functional and the shadow variable.
	Second, under the  representer assumption, we propose  nonparametric  estimation  that no longer involves parametrically  modelling the propensity score or  the outcome regression. Nonparametric estimation has been largely studied in the literature and has received broad acceptance with machine learning tools in many applications \citep{newey2003instrumental,chen2012estimation,kennedy2017non}.
	Because the solution to the representer equation may not be unique, we first   construct a consistent estimator of   the solution set. 
	We use the method of sieves to approximate unknown smooth functions as possible solutions and estimate corresponding coefficients by applying a minimum distance procedure, which has been routinely used in semiparametric and nonparametric econometric literature \citep{newey2003instrumental,ai2003efficient,santos2011instrumental,chen2012estimation}. 
	We then adapt the theory of extremum estimators to find from the estimated set a consistent estimator  of an appropriately chosen solution. Based on such an estimator, we propose a representer-based estimator for the  mean functional.
	Under certain regularity conditions, we establish  consistency and asymptotic normality for the proposed estimator.
	The proposed estimator is  shown to be   locally efficient for the mean functional under our shadow variable model.
	Besides,  due to the nonuniqueness of solutions to the representer equation,  
	the asymptotic results cannot be simply obtained by  applying  the standard  nonparametric sieve theory and  additional techniques are required.
	
	The remainder of this paper is organised as follows. 
	In Section \ref{sec:identification}, we provide a necessary and sufficient identifying condition for a generic mean functional within the shadow variable framework, and also introduce a representer assumption that is shown to be necessary for $\sqrt{n}$-estimability of the mean functional.
	In Section \ref{sec:estimation-inference},  we develop a model-free estimator for the mean functional,
	establish the asymptotic theory, discuss its semiparametric   efficiency, and provide computational details.
	In Section \ref{sec:numerical-studies}, we study the finite-sample performance of the proposed approach via both simulation studies and a real data example about  home pricing. 
	We conclude with a discussion in Section \ref{sec:discussion} and relegate proofs to the supporting information.
	
		\section{Identification}\label{sec:identification}
	
	Let $X$ denote  a vector of  fully observed covariates,   $Y$   the outcome
	variable that is subject to  missingness, and   $R$   the missingness
	indicator with $R=1$ if $Y$ is observed and $R=0$  otherwise.
	The missingness process may depend on the missing values.
	We let $f(\cdot)$ denote the probability density or mass function of a random variable (vector).
	The observed data contain $n$ independent and identically distributed  realizations of $(R,X,Y,Z)$ with the values of $Y$ missing for $R=0$.
	We are interested in identifying and making inference about the mean functional  $\mu = E\{\tau(X,Y)\}$ for a generic function $\tau(\cdot)$. In particular, when $\tau(X,Y)=Y$, the parameter $\mu=E(Y)$ corresponds to the population outcome mean that is of particular interest in sampling survey and causal inference. 	This setup also covers the   ordinary least squares regression problem  if we choose $\mu=E(X_j Y)$ for any component $X_j$ of $X$.
	Suppose we observe an additional  shadow variable $Z$ that meets   the following assumption.
	\begin{assumption}[Shadow variable]\label{ass:shadow}
		(i) $Z\indep R\mid (X,Y)$;~~ (ii) $Z\nindep Y\mid X$.
	\end{assumption}
	Assumption~\ref{ass:shadow}  reveals that the shadow variable does not affect the  missingness process  given the covariates and outcome, and  it  is associated with    the outcome given the covariates. 
	This assumption has been   used   for adjustment of selection bias in  sampling   surveys \citep{kott2014calibration} and in  missing data literature \citep{d2010new,wang2014instrumental,miao2016varieties}.
	Examples and extensive discussions about the assumption can be found in \cite{zahner1992children,ibrahim2001using,miao2016varieties,miao2019identification}.
	Under Assumption~\ref{ass:shadow}, we have	
	\begin{equation}\label{eqn:odds-representer}
		E\{\gamma(X,Y)\mid R=1,X,Z\}=\beta(X,Z),
	\end{equation}
	where
	\[\gamma(X,Y)=\frac{f(R=0\mid X,Y)}{f(R=1\mid X,Y)},\quad\text{and}\quad \beta(X,Z)=\frac{f(R=0\mid X,Z)}{f(R=1\mid X,Z)}.\]
	Without further assumptions such as completeness conditions, $\gamma(X,Y)$ is generally not identifiable from~\eqref{eqn:odds-representer}. Nevertheless, 	by noting that $\mu=E\{R\tau(X,Y)+R\tau(X,Y)\gamma(X,Y)\}$, the expectation functional $\mu$ can be identified under conditions weaker than those required for identification of $\gamma(X,Y)$ itself.
	The discussions below are based on operators on Hilbert spaces.
	Let $L_2(X,Y,Z)$ denote the set of real valued functions of $(X,Y,Z)$ that are square integrable with respect to the conditional distribution of $(X,Y,Z)$ given $R=1$. Let $T: L_2(X,Y)\longmapsto L_2(X,Z)$ be the linear operator given by $T(\xi)=E\{\xi(X,Y)\mid R=1,X,Z\}$. Its adjoint $T': L_2(X,Z)\longmapsto L_2(X,Y)$ is the  linear map $T'(\eta)=E\{\eta(X,Z)\mid R=1,X,Y\}$. The range and null space of $T$ is denoted by $\mathcal{R}(T)$ and $\mathcal{N}(T) $, respectively. The orthogonal complement of a set $\mathcal{A}$ is denoted by $\mathcal{A}^\bot$ and its closure in the norm topology is $\text{cl}(\mathcal{A})$.
	\begin{theorem}\label{thm:identification-shadow}
		Under Assumption~\ref{ass:shadow}, $\mu$ is identifiable if and only if $\tau(X,Y)\in \mathcal{N}(T)^\bot$.
	\end{theorem}

	Since the definition of the operator $T$ involves only observed data, the identifying condition could be justified in principle without extra model assumptions on the missing data distribution. 	In contrast to  previous approaches  \citep{d2010new,miao2019identification}     that have to     identify  the full data distribution under varying completeness conditions, 
	our identification strategy   allows for a larger class of models where only the parameter of interest is uniquely identified even though the full data law may not be. One of the commonly-used completeness conditions is that  for any square-integrable function $g$, $E\{g(X,Y)\mid R=1, X,Z\}=0$ almost surely if and only if $g(X,Y)=0$ almost surely. Under such circumstances,
	the function $\gamma(X,Y)$ defined after~\eqref{eqn:odds-representer} is identified, i.e., $\calN(T)=\{0\}$, then the functional $\mu$ is  identifiable because $\tau(X,Y)\in\mathcal{N}(T)^\bot=L_2(X,Y)$. 
	However, the identifying condition $\tau(X,Y)\in\mathcal{N}(T)^\bot$ does not suffice for 
	$\sqrt{n}$-estimability of $\mu$, particularly for functionals that lie in the boundary of $\mathcal{N}(T)^\bot$. In fact, following the proof of Lemma~4.1 in~\cite{severini2012efficiency}, we can show that
	Assumption~\ref{ass:representer} is necessary for $\sqrt{n}$-estimability of $\mu$
	under some regularity conditions.

	
	

	\begin{assumption}[Representer]\label{ass:representer}
		$\tau(X,Y)\in \calR(T')$, i.e.,
		there exists a function $\delta_0(X,Z) $ such that
		\begin{equation}\label{eqn:representer}
			E \big\{ \delta_0(X,Z) \mid R=1,X,Y \big\} =\tau(X,Y).  
		\end{equation}
	\end{assumption}
	Note that $\mathcal{N}(T)^\bot=\text{cl}\{\calR(T')\}$, the representer assumption   naturally strengthens the identifying condition $\tau(X,Y)\in\mathcal{N}(T)^\bot$. 
	Assumption \ref{ass:representer} is not only a sufficient  condition for identifiability of $\mu$, but also  necessary  for $\sqrt{n}$-estimability of $\mu$.
	Equation \eqref{eqn:representer}   relates the function $\tau(x,y)$ and the shadow variable  via the  representer function $\delta_0(x,z)$.
	The equation is  a Fredholm integral equation of the first kind, and the requirement for existence of solutions to \eqref{eqn:representer} is mild. Assumption~\ref{ass:representer} is nearly necessary for the completeness condition. Specifically, if the completeness condition holds, then  the solution to \eqref{eqn:representer} exists under the technical conditions~\ref{cond:A1}--\ref{cond:A3} given in the Appendix. Similar  conditions for existence of solutions to the Fredholm integral equation of the first kind are discussed in \cite{miao2018identify,carrasco2007,cui2020semiparametric}. 
	We give some specific examples that Assumption~\ref{ass:representer} holds. If  there exists  some transformation of $Z$ such that $E\{\lambda(Z)\mid X,Y\} = \alpha(X) + \beta(X) \tau(X,Y)$ and $\beta(x)\neq 0$, then Assumption \ref{ass:representer}  is met  with $\delta_0(X,Z) = \{\lambda(Z)-\alpha(X)\}/\beta(X)$.  
	As a special case, $\lambda(Z)=Z$ when $E(Z\mid X,Y)$ is linear in $\tau(X,Y)$. For simplicity,  we may drop the arguments in $\delta_0(X,Z)$ and directly
	use $\delta_0$ in what follows, and notation for other functions are treated in a similar way.
	
	Note that Assumption~\ref{ass:representer} only requires  the existence of  solutions to equation \eqref{eqn:representer}, but   not  uniqueness. 
	For instance, if both $Z$ and $Y$ are binary,  then $\delta_0$ is unique and 
	\begin{eqnarray*}
		\delta_0(X,Z)=\frac{Z\{\tau(X,1)-\tau(X,0)\}-f_{0}(X)\tau(X,1)+f_{1}(X)\tau(X,0)}{f_{1}(X)-f_{0}(X)},
	\end{eqnarray*}
	where $f_{0}(X)=f(Z=1\mid R=1,X,Y=0)$ and $f_{1}(X)=f(Z=1\mid R=1,X,Y=1)$.
	However, if $Z$ has more levels than $Y$, $\delta_0$ may not be unique.

	\begin{corollary}\label{cor:identifiable}
		Under Assumptions~\ref{ass:shadow} and \ref{ass:representer},  $\mu$ is  identifiable, and 
		\[
		\mu=E\big\{R\tau(X,Y) + (1-R)\delta_0(X,Z)\big\}.
		\]
	\end{corollary}
	
	From Corollary~\ref{cor:identifiable}, even if $\delta_0$ is not   uniquely determined, all solutions to Assumption \ref{ass:representer}  must result in an identical value of   $\mu$.
	Moreover,  this  identification result  does not require  identification of the full data  distribution  $f(R,X,Y,Z)$.
	In fact, identification of $f(R,X,Y,Z)$ is not ensured under   Assumptions \ref{ass:shadow} and \ref{ass:representer} only; see Example~\ref{exam:counter}.
	To our knowledge, the identifying Assumptions \ref{ass:shadow}--\ref{ass:representer} are so far the weakest for the shadow variable approach. Besides that, Assumption \ref{ass:representer} is also necessary for $\sqrt{n}$-estimability of $\mu$ within the shadow variable framework. 
	We further illustrate Assumption~\ref{ass:representer} with the following example.
	
	
	\begin{example}\label{exam:counter}
		Consider the following two models:
		
		Model 1: $Y\sim U(0,1)$, $Z\mid y\sim \Bern(y)$, and $f(R=1\mid y,z)=4y^2(1-y)$, where $\Bern(y)$ denotes Bernoulli distribution with probability $y$.
		
		Model 2: $Y\sim \Be(2,2)$, $Z\mid y\sim \Bern(y)$, and $f(R=1\mid y, z)=2y/3$, where $\Be(2,2)$ denotes Beta distribution with parameters 2 and 2.
		
		Suppose we are interested in estimating the outcome mean $\mu=E(Y)$. It is easy to verify that the above two models satisfy Assumption~\ref{ass:representer} by choosing $\delta_0(X,Z)=Z$. These two models imply the same outcome mean $E(Y)=1/2$ and the same observed data distribution, because $f(R=1,y,z)=4y^2(1-y)\{zy+(1-z)(1-y)\}$ and $f(z)=1/2$ in these two models.
		However, the full data distributions of these two models are different.
	\end{example}

	\section{Estimation, inference, and computation}\label{sec:estimation-inference}

In this section, we provide a novel estimation procedure without modeling the propensity score or outcome regression.  Previous approaches often require fully or partially parametric models for at least one of them. For example, \cite{qin2002estimation} and \cite{wang2014instrumental} assumed a fully parametric model for the propensity score;
\cite{kim2011semiparametric} and \cite{shao2016semiparametric} relaxed their assumption and considered a semiparametric exponential tilting model for the propensity; \cite{miao2016varieties} proposed doubly robust estimation methods by either requiring a parametric propensity score or an outcome regression to be correctly specified. Our approach aims to be more robust than existing methods by avoiding (i) point identification of the full data law under more stringent conditions, and (ii) over-reliance on parametric assumptions either for identification or for estimation.

As implied by Corollary~\ref{cor:identifiable}, any solution to~\eqref{eqn:representer} provides a valid $\delta_0$ for recovering the parameter $\mu$.
Suppose that all such solutions belong to a set $\Delta$ of smooth functions, with specific requirements for smooth functions given in Definition~\ref{def:smooth}. Then the set of solutions to \eqref{eqn:representer} is denoted by
\begin{equation}\label{eqn:identifiable-set}
	\Delta_0=\Big\{\delta\in \Delta: E \big\{ \delta(X,Z) \mid R=1,X,Y \big\} =\tau(X,Y) \Big\}.
\end{equation}

For estimation and inference about $\mu$, we need to construct a  consistent estimator for some fixed $\delta_0\in\Delta_0$. 
If $\Delta_0$ were known, then we would simply select one element $\delta_0$  from the set and use this element to estimate $\mu$. Unfortunately, the solution set $\Delta_0$ is unknown, and
the lack of identification of $\delta_0$ presents important technical challenges. Directly solving~\eqref{eqn:representer} does not generally yields a consistent estimator for some fixed $\delta_0$. Instead, by noting that the solution set $\Delta_0$ is identified, we aim to obtain an estimator $\widehat\delta_0$ in the following two steps: 
first, construct a consistent estimator $\widehat\Delta_0$ for the set $\Delta_0$; second, carefully select $\widehat\delta_0\in\widehat\Delta_0$ such that it is a consistent estimator for a fixed element $\delta_0\in\Delta_0$. 

\subsection{Estimation of the solution set $\Delta_0$}

Define the criterion function
\begin{equation*}
	Q(\delta)=E\Big[R\big\{E(\tau(X,Y)-\delta(X,Z)\mid R=1,X,Y) \big\}^2 \Big].
\end{equation*}
Then the solution set $\Delta_0$ in~\eqref{eqn:identifiable-set} is   equal to the set of zeros of  	$Q(\delta)$, i.e.,
\begin{equation*}
	\Delta_0=\big\{\delta\in \Delta: Q(\delta)=0 \big\},
\end{equation*}
and hence, estimation of $\Delta_0$ is equivalent to estimation of zeros of $Q(\delta)$. 
This can be accomplished  with the approximate minimizers of a sample analogue of $Q(\delta)$ \citep{chernozhukov2007estimation}. 

We adopt a method of sieves approach to construct a sample analogue function $Q_n(\delta)$ for $Q(\delta)$ and a corresponding approximation $\Delta_n$ for  $\Delta$.
Let $\{\psi_q(x,z)\}_{q=1}^{\infty}$ denote a sequence of known approximating functions of $x$ and $z$, and 
\begin{equation}\label{eqn:Deltan-sieve}
	\Delta_n=\bigg\{\delta\in \Delta: \delta(x,z)=\sum_{q=1}^{q_n}\beta_q \psi_q(x,z) \bigg\}
\end{equation}
for some known $q_n$ and unknown parameters $\{\beta_q\}_{q=1}^{q_n}$. 
The construction of $Q_n$ entails a nonparametric estimator of  conditional expectations. 
Let  $\{\phi_k(x,y)\}_{k=1}^{\infty}$ be a sequence of known approximating functions of $x$ and $y$.  
Denote the vector of the first $k_n$ terms of the basis functions by
\begin{equation*}
	{\bphi}(x,y) = \big\{\phi_1(x,y),\ldots,\phi_{k_n}(x,y)\big\}^\T,
\end{equation*}
and let
\begin{equation*}
	\bPhi=\big\{\bphi(X_1,Y_1),\ldots,\bphi(X_n,Y_n)\big\}^{\T},\qquad \bLambda= \text{diag}\big(R_1,\ldots,R_n\big).
\end{equation*}
For a generic random variable $B=B(X,Y,Z)$ with realizations $\{B_i=B(X_i,Y_i,Z_i) \}_{i=1}^n$, 
the nonparametric sieve estimator  of $E(B\mid R=1, x,y)$ is  obtained by the linear regression of $B$ on the vector $\bphi(X,Y)$ with  observed data, i.e.,
\begin{equation}\label{eqn:sieve-conditional-expectation}
	\widehat E(B\mid R=1, X,Y) = \bphi^\T(X,Y)(\bPhi^\T\bLambda\bPhi)^{-1}\sum_{i=1}^n R_i\bphi(X_i,Y_i) B_i.
\end{equation}
Then the sample analogue $Q_n$  of $Q$ is  
\begin{equation}\label{eqn:sample-analogue-Q}
	\begin{aligned}
		Q_n(\delta)&=\frac{1}{n}\sum_{i=1}^n R_i  \widehat e^2(X_i,Y_i,\delta),
	\end{aligned}
\end{equation}
with 
\begin{equation}\label{eqn:hat-e}
	\begin{aligned}
		\widehat e(X_i,Y_i,\delta)&=\widehat E\big\{\tau(X,Y)-\delta(X,Z)\mid R=1,X_i,Y_i\big\},
	\end{aligned}
\end{equation}
where the explicit expression of $\widehat e(X_i,Y_i,\delta)$ is obtained from \eqref{eqn:sieve-conditional-expectation}.
Finally, the proposed estimator of $\Delta_0$ is  
\begin{equation}\label{eqn:estimator-for-Delta0}
	\widehat\Delta_0=\big\{\delta\in \Delta_n: Q_n(\delta)\leq c_n \big\},
\end{equation}
where $\Delta_n$ and $Q_n(\delta)$ are given in~\eqref{eqn:Deltan-sieve} and~\eqref{eqn:sample-analogue-Q}, respectively, and
$\{c_n\}_{n=1}^\infty$ is a sequence of small positive numbers converging to zero at an appropriate rate. The requirement on the rate of $c_n$ will be discussed later for theoretical analysis.

\subsection{Set consistency}
We establish the set consistency of $\widehat\Delta_0$ for $\Delta_0$ in terms of  Hausdorff distances. 
For a given norm $\Vert\cdot\Vert$,  the Hausdorff distance between two   sets  $\Delta_1,\Delta_2\subseteq \Delta$ is
\begin{equation*}
	\begin{aligned}
		d_H(\Delta_1,\Delta_2,\Vert \cdot\Vert) &= \max\big\{d(\Delta_1,\Delta_2),d(\Delta_2,\Delta_1)\big\},
	\end{aligned}
\end{equation*}
where $d(\Delta_1,\Delta_2)= \sup_{\delta_1\in\Delta_1}\inf_{\delta_2\in\Delta_2} \Vert \delta_1-\delta_2\Vert$ and $d(\Delta_2,\Delta_1)$ is defined analogously. 
Thus, $\widehat\Delta_0$ is consistent under the Hausdorff distance if both the maximal approximation error of $\widehat\Delta_0$ by $\Delta_0$ and of $\Delta_0$ by $\widehat\Delta_0$ converge to zero in probability.

We consider two different norms for the Hausdorff distance: 
the pseudo-norm $\Vert\cdot\Vert_w$ defined by 
\begin{equation*}
	\Vert \delta\Vert_w^2 = E\Big[R\big\{ E(\delta(X,Z)\mid R=1,X,Y)\big\}^2 \Big],
\end{equation*}
and  the supremum norm $\Vert\cdot\Vert_{\infty}$ defined by
\[\Vert\delta\Vert_{\infty}=\sup_{x,z}|\delta(x,z)|.\]

%
%

From   the representer equation~\eqref{eqn:representer}, we have that for any $\widehat\delta_0\in\widehat\Delta_0$ and $\delta_0,\delta\in\Delta_0$, $ \Vert \widehat\delta_0-\delta_0\Vert_w=\Vert\widehat\delta_0-\delta\Vert_w$. Hence,
\begin{equation}\label{eqn:element-set-relation}
	\Vert \widehat\delta_0-\delta_0\Vert_w=\inf_{\delta\in\Delta_0}\Vert\widehat\delta_0-\delta\Vert_w\leq d_H(\widehat\Delta_0,\Delta_0,\Vert\cdot\Vert_w).
\end{equation}
This result  implies that we can obtain the convergence rate of  $\Vert\widehat\delta_0-\delta_0\Vert_w$ by deriving that of $d(\widehat\Delta_0,\Delta_0,\Vert\cdot\Vert_w)$. 
However,   the identified set $\Delta_0$ is an equivalence class under the pseudo-norm, and the convergence under  $d_H(\widehat\Delta_0,\Delta_0,\Vert\cdot\Vert_w)$ does not suffice to consistently estimate a given element $\delta_0\in\Delta_0$. 
Whereas the supremum norm $\Vert\cdot\Vert_{\infty}$  is able to differentiate between elements in $\Delta_0$,  
and   $d_H(\widehat\Delta_0,\Delta_0,\Vert\cdot\Vert_{\infty})=o_p(1)$ under certain regularity condition as we will show later.

We make the following assumptions to guarantee that $\widehat\Delta_0$ is consistent under the metric $d_H(\widehat\Delta_0,\Delta_0,\Vert\cdot\Vert_\infty)$ and to obtain the rate of convergence for $\widehat\Delta_0$ under the weaker metric $d_H(\widehat\Delta_0,\Delta_0,\Vert\cdot\Vert_w)$.

\begin{assumption}\label{ass:data-distribution}
	The vector of covariates $X\in \R^{d}$ has support $[0,1]^{d}$, and the outcome $Y\in\R$ and the shadow variable $Z\in \R$ have compact supports.
\end{assumption}

Assumption~\ref{ass:data-distribution} requires  $(X,Y,Z)$ to have compact supports, and  without loss of generality, 
we assume that $X$ has been  transformed such that the support  is $[0,1]^{d}$. 
These are standard conditions that are usually required in the semiparametric literature.
Although $Y$ and $Z$ are also required to have compact support, the proposed approach may still be applicable if the supports are infinite with sufficiently thin tails.
For instance,  in our simulation studies where the variables $Y$ and $Z$ are drawn from a normal distribution in Section~\ref{sec:numerical-studies},  the proposed approach continues to perform quite well.

We next impose  restrictions on the smoothness of functions in the set $\Delta$. We use  the following Sobolev norm to   characterize the smoothness of functions.
\begin{definition}\label{def:smooth}
	For a generic function  $\rho(w)$ defined on $w\in\R^d$,  we define
	\begin{equation*}
		\Vert\rho\Vert_{\infty,\alpha}=\max_{|\lambda|\leq \underline{\alpha}}\sup_w|D^{\lambda}\rho(w)|+\max_{\lambda=\underline{\alpha}}\sup_{w\neq w'}
		\frac{D^\lambda\rho(w)-D^{\lambda}\rho(w')}{\Vert w-w'\Vert^{\alpha-\underline{\alpha}}},
	\end{equation*}
	where $\lambda$ be a $d$-dimensional vector of nonnegative integers,  $|\lambda|=\sum_{i=1}^{d}\lambda_i$,
	$\underline{\alpha}$ denotes the largest integer smaller than $\alpha$,
	$D^{\lambda}\rho(w)=\partial^{|\lambda|}\rho(w)/\partial w_1^{\lambda_1}\ldots \partial w_{d}^{\lambda_{d}}$, and  $D^0\rho(w)=\rho(w)$.
\end{definition}
A function $\rho$ with $\Vert\rho\Vert_{\infty,\alpha}<\infty$ has uniformly bounded partial derivatives up to order $\underline{\alpha}$; besides, the $\underline{\alpha}$th partial derivative of this function is Lipschitz of order $\alpha-\underline{\alpha}$. 

\begin{assumption}\label{ass:parameter-space}
	The following conditions hold:
	\begin{itemize}
		\item[(i)]  
		$\sup_{\delta\in\Delta} \Vert\delta\Vert_{\infty,\alpha}<\infty$ for some $\alpha>(d+1)/2$;   in addition, $\Delta_0\neq\emptyset$,  and both $\Delta_n$ and $\Delta$ are closed;
		\item[(ii)] for every $\delta\in\Delta$, there is $\Pi_n\delta\in\Delta_n$ such that $\sup_{\delta\in\Delta}\Vert\delta-\Pi_n\delta\Vert_{\infty}=O(\eta_n)$ for some $\eta_n=o(1)$.
	\end{itemize}
\end{assumption}

Assumption~\ref{ass:parameter-space}(i) requires that each function $\delta\in\Delta$ is sufficiently smooth and bounded. The closedness condition in this assumption and Assumption~\ref{ass:data-distribution} together imply that $\Delta$ is compact under $\Vert\cdot\Vert_{\infty}$.
It is well known that solving integral equations as in~\eqref{eqn:representer} is an ill-posed inverse problem. The ill-posedness  due to noncontinuity of the solution and difficulty of computation can have a severe impact on the consistency and convergence rates of estimators. 
The compactness condition   is imposed to ensure that the consistency of the proposed estimator under $\Vert\cdot\Vert_{\infty}$ is not affected by
the ill-posedness.
Such a compactness condition is commonly made in the nonparametric and semiparametric literature; see, e.g., \cite{newey2003instrumental}, \cite{ai2003efficient}, and \cite{chen2012estimation}. 
Alternatively, it is possible to address the ill-posed problem by employing a regularization approach as in \cite{horowitz2009semiparametric} and \cite{darolles2011nonparametric}.

Assumption~\ref{ass:parameter-space}(ii) quantifies the approximation error of functions in $\Delta$ by the sieve space $\Delta_n$. This condition is satisfied by many commonly-used function spaces  (e.g., H\"{o}lder space), whose elements are sufficiently smooth, and by popular sieves (e.g., power series, splines). 
For example, consider the function set $\Delta$ with  $\sup_{\delta\in\Delta} \Vert\delta\Vert_{\infty,\alpha}<\infty$.
If the sieve functions $\{\psi_q(x,z)\}_{q=1}^{\infty}$ are polynomials
or tensor product univariate splines, then uniformly on $\delta\in\Delta$, the
approximation error of $\delta$
by functions of the form $\sum_{q=1}^{q_n}\beta_q\psi_q(x,z)\in\Delta_n$ under $\Vert\cdot\Vert_{\infty}$
is of the order $O\{q_n^{-\alpha/(d+1)}\}$.
Thus,  Assumption~\ref{ass:parameter-space}(ii) is met  with $\eta_n=q_n^{-\alpha/(d+1)}$; see \cite{chen2007large} for further discussion.

\begin{assumption}\label{ass:expectation-basis-function}
	The following conditions hold:
	\begin{itemize}
		\item[(i)] the smallest and largest eigenvalues of $E\{R\bphi(X,Y) \bphi(X,Y)^\T\}$ are bounded above and  away from zero for all $k_n$;
		\item[(ii)] for every $\delta\in\Delta$, there is a $\bm\pi_n(\delta)\in\R^{k_n}$ such that
		\[
		\sup_{\delta\in\Delta}\Vert E\{\delta(X,Z)\mid r=1,x,y\}-\bphi^\T(x,y){\bm\pi}_n(\delta)\Vert_{\infty}=O\Big(k_n^{-\frac{\alpha}{d+1}}\Big);
		\]
		\item[(iii)]  $\xi_n^2 k_n=o(n)$, where $\xi_n = \sup_{x,y}\Vert\bphi(x,y)\Vert_2$.
		
	\end{itemize}
\end{assumption}

Assumption~\ref{ass:expectation-basis-function} bounds the second  moment matrix of the approximating functions away from singularity, presents a uniform approximation error of the series estimator to the conditional mean function, and restricts the magnitude of the series terms. These conditions are standard  for  series estimation of conditional mean functions; see, e.g., \cite{newey1997convergence}, \cite{ai2003efficient}, and \cite{huang2003local}. Primitive conditions are discussed below so that the rate requirements in this assumption hold.
Consider any $\delta\in\Delta$ satisfying Assumption~\ref{ass:parameter-space}, i.e., $\sup_{\delta\in\Delta} \Vert\delta\Vert_{\infty,\alpha}<\infty$. If the partial derivatives of $f(z\mid r=1,x,y)$ with respect to $(x,y)$ are continuously differentiable up to order $\underline{\alpha}+1$, then under Assumption~\ref{ass:data-distribution}, we have $\sup_{\delta}
\Vert E\{\delta(X,Z)\mid R=1,x,y\}\Vert_{\infty,\alpha}<\infty$. In addition, if the sieve functions $\{\phi_k(x,y)\}_{k=1}^{k_n}$ are polynomials or tensor product univariate splines, then
by similar arguments after Assumption~\ref{ass:parameter-space}, we conclude that
the approximation error under $\Vert\cdot\Vert_{\infty}$ is of the order $O\{k_n^{-\alpha/(d+1)}\}$ uniformly on $\delta\in\Delta$. Verifying Assumption~\ref{ass:expectation-basis-function}(iii) depends on the relationship between $\xi_{n}$ and $k_n$. For example, if $\{\phi_k(x,y)\}_{k=1}^{k_n}$ are tensor product univariate splines, then  $\xi_{n}=O\{ k_n^{(d+1)/2}\}$.


Write $c_n$ in \eqref{eqn:estimator-for-Delta0} by $b_n/a_n$ with appropriate sequences $a_n$ and $b_n$, and define $\lambda_n=k_n/n+k_n^{-2\alpha/(d+1)}+\eta_n^2$.
\begin{theorem}\label{thm:set-convergence}
	Suppose that Assumptions~\ref{ass:data-distribution}--\ref{ass:expectation-basis-function} hold.
	If $a_n=O(\lambda_n^{-1})$, $b_n\rightarrow\infty$ and $b_n=o(a_n)$, 
	Then 
	\begin{equation*}
		d_H\big(\widehat\Delta_0,\Delta_0,\Vert\cdot\Vert_{\infty}\big)=o_p(1),\quad \text{and}\quad
		d_H\big(\widehat\Delta_0,\Delta_0,\Vert\cdot\Vert_{w}\big)=O_p\big(c_n^{1/2}\big).
	\end{equation*}
\end{theorem}

Theorem~\ref{thm:set-convergence} shows the consistency of $\widehat\Delta_0$ under the supremum-norm metric $d_H(\widehat\Delta_0,\Delta_0,\Vert\cdot\Vert_{\infty})$ and establishes
the rate of convergence of $\widehat\Delta_0$ under the weaker pseudo-norm metric  $d_H(\widehat\Delta_0,\Delta_0,\Vert\cdot\Vert_{w})$.
In particular, if we let $k_n^3=o(n)$, $k_n^{-3\alpha/(d+1)}=o(n^{-1})$, and $\eta_n=o(n^{-1/3})$ as imposed in Assumption~\ref{ass:ASN-rate} in the next subsection, then $\lambda_n=o(n^{-2/3})$ or $\lambda_n^{-1}n^{-2/3}\rightarrow\infty$. We  take $a_n=\lambda_n^{-1/2}n^{1/3}\rightarrow\infty$ and $b_n=a_n^{1/2}/n^{1/3}$. Thus, $a_n=\lambda_n^{-1}(\lambda_nn^{2/3})^{1/2}=o(\lambda_n^{-1})$, $b_n=(\lambda_n^{-1}n^{-2/3})^{1/4}\rightarrow\infty$, and $b_n=a_n\cdot a_n^{-1/2}n^{-1/3}=o(a_n)$. In fact, under such rate requirements, we 
have
$n^{2/3}b_n=o(a_n)$ and $d_H(\widehat\Delta_0,\Delta_0,\Vert\cdot\Vert_{w})=o_p(n^{-1/4})$, which are sufficient 
to establish the asymptotic normality of the proposed estimator given in subsection~\ref{subsec:repa}.

\subsection{A representer-based estimator}\label{subsec:rep}

After we have obtained a consistent estimator $\widehat\Delta_0$ for $\Delta_0$,  we remain to  select an estimator from $\widehat\Delta_0$  such that it converges to a unique element belonging to $\Delta_0$. 
We adapt the theory of extremum estimators to achieve this goal.
Let $M:\Delta\rightarrow \R$ be a population criterion functional that attains a unique minimum $\delta_0$ on   $\Delta_0$ and 
$M_n(\delta)$ be its sample analogue.
We then choose  the minimizer of $M_n(\delta)$ over the estimated solution set $\widehat\Delta_0$, denoted by 
\begin{equation}\label{eqn:hat-delta0}
	\widehat\delta_0\in\argmin_{\delta\in\widehat\Delta_0}M_n(\delta),
\end{equation} 
which   is  expected to converge to  the unique minimum $\delta_0$ of $M(\delta)$ on $\Delta_0$. 
\begin{assumption}\label{ass:extremum}
	The function set $\Delta$ is convex;
	the functional $M:\Delta\rightarrow \R$ is strictly convex and attains a unique minimum at $\delta_0$ on $\Delta_0$;
	its sample analogue $M_n:\Delta\rightarrow \R$ is continuous and  $\sup_{\delta\in\Delta}|M_n(\delta)-M(\delta)|=o_p(1)$. 
\end{assumption}

One example  of particular interest is
\begin{equation*}
	M(\delta)=E\Big[\big\{(1-R)\delta(X,Z)\big\}^2\Big].
\end{equation*}
This is a convex functional with respect to $\delta$. 
In addition, since $E\{(1-R)\delta_0(X,Z)\}=E\{(1-R)\tau(X,Y)\}$ for any $\delta_0\in\Delta_0$, the minimizer of $M(\delta)$ on $\Delta_0$ in fact minimizes the variance of $(1-R)\delta_0(X,Z)$ among $\delta_0\in\Delta_0$.
Its sample analogue   is
\begin{equation*}\label{eqn:Mn-function}
	M_n(\delta)=\frac{1}{n}\sum_{i=1}^n (1-R_i)\delta^2(X_i,Z_i). 
\end{equation*}
Under Assumptions~\ref{ass:data-distribution}--\ref{ass:parameter-space}, one can show that the function class $\{(1-R)\delta:\delta\in\Delta\}$  is a Glivenko-Cantelli class, and thus $\sup_{\delta\in\Delta}|M_n(\delta)-M(\delta)|=o_p(1)$.

\begin{theorem}\label{thm:element-rate}
	Suppose that Assumptions~\ref{ass:data-distribution}--\ref{ass:extremum} hold.  Then
	\begin{equation*}
		\Vert\widehat\delta_0-\delta_0\Vert_{\infty}=o_p(1),
	\end{equation*}
	where $\widehat\delta_0$ is defined through~\eqref{eqn:hat-delta0} and $\delta_0$ is defined in Assumption~\ref{ass:extremum}.
	In addition,  if $a_n=O(\lambda_n^{-1})$, $b_n\rightarrow\infty$ and $b_n=o(a_n)$, we then have
	%
	\begin{equation*}
		\Vert\widehat\delta_0-\delta_0\Vert_w=O_p(c_n^{1/2}).
	\end{equation*}
\end{theorem}

Theorem~\ref{thm:element-rate} implies that by choosing an appropriate  function $M(\delta)$, 
it is possible to construct a consistent estimator $\widehat\delta_0$ for some unique element $\delta_0\in\Delta_0$ in terms of supremum norm $\Vert\cdot\Vert_{\infty}$ and further obtain its rate of convergence under the weaker pseudo-norm $\Vert\cdot\Vert_w$.

Based on the   estimator $\widehat\delta_0$ given in~\eqref{eqn:hat-delta0}, we obtain the following representer-based estimator $\widehat\mu_{\rep}$ of $\mu$:
\begin{equation}\label{eqn:representer-based-estimator}
	\widehat\mu_{\rep}=\frac{1}{n}\sum_{i=1}^n \Big\{R_i\tau(X_i,Y_i) +(1-R_i)\widehat\delta_0(X_i,Z_i)\Big\}.
\end{equation}
Below we discuss the asymptotic expansion of the estimator $\widehat\mu_{\rep}$.

Let $\overline \Delta$   be the closure of the linear span of $\Delta$ under $\Vert\cdot\Vert_w$, which is a Hilbert space with inner product: 
\begin{equation*}
	\langle\delta_1,\delta_2\rangle_w=E\Big[ RE\big\{\delta_1(X,Z)\mid R=1,X,Y\big\}E\big\{\delta_2(X,Z)\mid R=1,X,Y\big\} \Big]
\end{equation*}
for any $\delta_1,\delta_2\in\overline \Delta$.


\begin{assumption}\label{ass:ASN-rate}
	The following conditions hold:
	\begin{itemize}
		\item[(i)] there exists a function $  h_0\in \Delta$ such that  
		\begin{equation*}\label{eqn:tilde-H}
			\langle h_0, \delta\rangle_w=E\big\{(1-R)\delta(X,Z)\big\}~\text{for all $\delta \in\overline \Delta$}.
		\end{equation*}
		
		\item[(ii)] $\eta_n=o(n^{-1/3})$,
		$k_n^{-3\alpha/(d+1)}=o(n^{-1})$, $k_n^3=o(n)$,  $\xi_{n}^2k_n^2=o(n)$, and $\xi_{n}^2k_n^{-2\alpha/(d+1)}=o(1)$.
	\end{itemize}
\end{assumption}

Note that the linear functional $\delta\longmapsto E\{(1-R)\delta(X,Z)\}$ is continuous under $\Vert\cdot\Vert_w$. Hence, by the Riesz representation theorem, there exists a unique $  h_0\in\overline \Delta$ (up to an equivalence class in $\Vert\cdot\Vert_w$)
such that  $\langle  h_0,\delta\rangle_w=E\{(1-R)\delta(X,Z)\}$ for all $\delta\in\overline \Delta$.
However, Assumption~\ref{ass:ASN-rate}(i) further requires that this equivalence class must contain  at least one  element that falls in    $\Delta$.
A primitive condition for  Assumption~\ref{ass:ASN-rate}(i) is that  the inverse probability weight also has a smooth representer:  if
\begin{equation}\label{eqn:second-representer-equation}
	E\big\{    h_0(X,Z) + 1\mid R=1,X,Y\big\}=\frac{1}{f(R=1\mid X,Y)},
\end{equation}
then $h_0$ satisfies  Assumption~\ref{ass:ASN-rate}(i).

Assumption~\ref{ass:ASN-rate}(ii) imposes some rate requirements, which can be satisfied as long as the function classes being approximated in Assumptions~\ref{ass:parameter-space} and \ref{ass:expectation-basis-function} are sufficiently smooth.

\begin{theorem}\label{thm:asymptotic-expansion}
	Suppose that Assumptions~\ref{ass:data-distribution}--\ref{ass:ASN-rate} hold.  
	We have  that
	\begin{equation*}
		\begin{aligned}
			\sqrt{n}(\widehat\mu_{\rep}-\mu)=&\frac{1}{\sqrt{n}}\sum_{i=1}^n\Big[(1-R_i)\delta_0(X_i,Z_i)+R_i\tau(X_i,Y_i)+
			R_iE\big\{h_0(X,Z)\mid R=1,X_i,Y_i\big\} \\
			&	~~~~~~~~~~~~		\times\big\{\tau(X_i,Y_i)-\delta_0(X_i,Z_i) \big\} -\mu\Big]-\sqrt{n}r_n(\widehat\delta_0)+o_p(1),
		\end{aligned}
	\end{equation*}
	with
	\begin{equation}\label{eq:bias}
		r_n(\widehat\delta_0)=\frac{1}{n}\sum_{i=1}^nR_i\widehat E\big\{\Pi_n h_0(X,Z)\mid R=1,X_i,Y_i\big\}\widehat e(X_i,Y_i,\widehat\delta_0),
	\end{equation}
	where  $\Pi_nh_0\in\Delta_n$ approximates $h_0$ as given in Assumption~\ref{ass:parameter-space}(ii), $\widehat E(\cdot)$ and $\widehat e(\cdot)$ are defined in~\eqref{eqn:sieve-conditional-expectation} and~\eqref{eqn:hat-e}, respectively.
\end{theorem}
Theorem~\ref{thm:asymptotic-expansion} reveals an asymptotic expansion of $\widehat\mu_{\rep}$. 
However, the estimator $\widehat\mu_{\rep}$ is not necessarily  asymptotically normal as the bias term $\sqrt{n}r_n(\widehat\delta_0)$ may not be asymptotically negligible. 
In the next subsection, we propose a debiased estimator which is regular and asymptotically normal. We further establish that the debiased estimator is semiparametric locally efficient under a shadow variable model at a given submodel where a key completness condition holds.

\subsection{A debiased  semiparametric locally efficient estimator}\label{subsec:repa}
Note that only $\Pi_nh_0$ is unknown in the bias term $r_n(\widehat\delta_0)$   in \eqref{eq:bias}.
We propose to construct an estimator of $\Pi_nh_0$ and then subtract the bias to obtain an  estimator of $\mu$ that is asymptotically normal. 
We  define the criterion function: 
\begin{equation*}
	C(\delta)=E\Big[R\big\{E(\delta(X,Z)\mid R=1,X,Y)\big\}^2 \Big]-2E\big\{(1-R)\delta(X,Z)\big\}, \quad\delta\in\Delta
\end{equation*}
and its sample analogue,
\begin{equation*}\label{eqn:Cn}
	\begin{aligned}
		C_n(\delta)=\frac{1}{n}\sum_{i=1}^n R_i\Big[\widehat E\{\delta(X,Z)\mid R=1,X_i,Y_i \} \Big]^2-\frac{2}{n}\sum_{i=1}^n(1-R_i)\delta(X_i,Z_i),\quad\delta\in\Delta.
	\end{aligned}
\end{equation*}
Since $E\{(1-R)\delta(X,Z)\}=\langle h_0,\delta\rangle_w$ by Assumption~\ref{ass:ASN-rate}, it follows that $C(\delta)=\Vert\delta- h_0\Vert_w^2-\Vert h_0\Vert_w^2$.  Thus, $ h_0$ is the unique minimizer of $\delta \longmapsto C(\delta)$ up to the equivalence class in $\Vert\cdot\Vert_w$. In addition, since $ h_0$ and $\Pi_n h_0$ are close under the metric $\Vert\cdot\Vert_{\infty}$ by Assumption~\ref{ass:parameter-space}(ii), we then
define the estimator for $\Pi_n h_0$ by:
\begin{equation}\label{eqn:hat-h}
	\begin{aligned}
		\widehat h&\in\argmin_{\delta\in\Delta_n} C_n(\delta),
	\end{aligned}
\end{equation}
Given the estimator $\widehat h$, the approximation to the bias term $r_n(\widehat\delta_0)$ is
\begin{equation}\label{eqn:estimated-remainder-term}
	\widehat r_n(\widehat\delta_0)=\frac{1}{n}\sum_{i=1}^n R_i\widehat E\Big\{\widehat h(X,Z)\mid R=1, X_i,Y_i\Big\}\widehat e(X_i,Y_i,\widehat\delta_0).
\end{equation}
\begin{lemma}\label{lem:remainder-estimator-rate}
	Suppose that Assumptions~\ref{ass:data-distribution}--\ref{ass:expectation-basis-function} and \ref{ass:ASN-rate} hold. Then it follows that
	\begin{equation*}
		\sup_{\widehat\delta_0\in\widehat\Delta_0}\big|\widehat r_n(\widehat\delta_0)-r_n(\widehat\delta_0)\big|=O_p\Bigg[c_n^{1/2}\bigg\{\bigg(\frac{k_n}{n} \bigg)^{1/4}+k_n^{-\frac{\alpha}{2(d+1)} } \bigg\}\Bigg].
	\end{equation*}
\end{lemma}
This lemma establishes the rate of convergence of $\widehat r_n(\widehat\delta_0)$ to $r_n(\widehat\delta_0)$ uniformly on $\widehat\Delta_0$.
If $c_n$ converges to zero sufficiently fast, then $\sup_{\widehat\delta_0\in\widehat\Delta_0}\sqrt{n}|\widehat r_n(\widehat\delta_0)-r_n(\widehat\delta_0)|=o_p(1)$. The rate conditions imposed in Assumption~\ref{ass:ASN-rate}(ii) guarantee that such a choice of $c_n$ is   feasible. 
As a result, Theorem~\ref{ass:ASN-rate} and Lemma~\ref{lem:remainder-estimator-rate} imply that it is possible to construct a debiased estimator that is $\sqrt{n}$-consistent and asymptotically normal by subtracting the  estimated bias   $ \widehat r_n(\widehat\delta_0)$ from $\widehat\mu_{\rep}$:
\begin{equation}\label{eqn:representer-based-estimator-adjustment}
	\widehat\mu_{\repa}=\widehat\mu_{\rep}+\widehat r_n(\widehat\delta_0).
\end{equation}

\begin{theorem}\label{thm:ASN}
	Suppose that Assumptions~\ref{ass:data-distribution}--\ref{ass:ASN-rate} hold.  
	If $a_n=O(\lambda_n^{-1})$, $b_n\rightarrow\infty$ and $n^{2/3}b_n=o(a_n)$,
	then 
	$ \sqrt{n}(\widehat\mu_{\repa}-\mu)$ converges in distribution to $N(0,\sigma^2)$, where $\sigma^2 $ is the variance of 
	\begin{equation}\label{eqn:influence-function}
		(1-R)\delta_0(X,Z)+R\tau(X,Y)+RE\{  h_0(X,Z)\mid R=1,X,Y\}\{\tau(X,Y)-\delta_0(X,Z)\}-\mu.
	\end{equation}
\end{theorem}
Based on~\eqref{eqn:influence-function}, one can easily obtain a consistent estimator of the asymptotic variance
\begin{equation*}
	\begin{aligned}
		\widehat\sigma^2 =& \frac{1}{n}\sum_{i=1}^n\bigg [(1-R_i)\widehat\delta_0^2(X_i,Z_i)+R_i\tau^2(X_i,Y_i)-\bigg\{\frac{1}{n}\sum_{i=1}^n(1-R_i)\widehat\delta_0(X_i,Z_i)+R_i\tau(X_i,Y_i) \bigg\}^2\\
		&~~~~~~~~~~~~
		+R_i\Big\{\widehat E(\widehat
		h(X,Z)\mid R=1,X_i,Y_i) \Big\}^2\Big\{\tau(X_i,Y_i)-\widehat\delta_0(X_i,Z_i) \Big\}^2 \bigg].
	\end{aligned}
\end{equation*}
Then given 
$\alpha\in(0,1)$, an asymptotic $100(1-\alpha)\%$ confidence interval is $[\widehat\mu_{\repa}-z_{\alpha}\widehat{\sigma}/\sqrt{n},\widehat\mu_{\repa}+z_{\alpha}\widehat{\sigma}/\sqrt{n}]$, where $z_{\alpha}=\Phi^{-1}(1-\alpha/2)$.
The formula \eqref{eqn:influence-function} presents the influence function for $\widehat\mu_{\repa}$. 
The influence function  is locally efficient in the sense that  it attains the semiparametric efficiency bound for the outcome mean under certain conditions in the semiparametric model $\mathcal M_{\np}$ defined through~\eqref{eqn:odds-representer}.
%

\begin{assumption}\label{ass:efficiency-bound}
	The following conditions hold:
	\begin{itemize}
		\item[(i)] 
		Completeness: (1) for any square-integrable function $\xi(x,y)$, $E\{\xi(X,Y)\mid R=1, X,Z\}=0$ almost surely if and only if $\xi(X,Y)=0$ almost surely; (2) for any square-integrable function $\eta(x,z)$, $E\{\eta(X,Z)\mid R=1, X,Y\}=0$ almost surely if and only if $\eta(X,Z)=0$ almost surely.
		\item[(ii)]
		Denote $\Omega(x,z)=E[\{\gamma(X,Y)-\beta(X,Z)\}^2\mid 
		R=1,X=x,Z=z]$. Suppose that $0<\inf_{x,z}\Omega(x,z)\leq \sup_{x,z}\Omega(x,z)<\infty$.
		\item[(iii)]
		The operator $T$ and its adjoint $T'$ defined through conditional expectations in Section~\ref{sec:identification} are both bounded.
	\end{itemize}
\end{assumption}

Under the completeness condition in Assumption~\ref{ass:efficiency-bound}(i), $\gamma(X,Y)$ is identifiable, $\delta_0(X,Z)$ and $h_0(X,Z)$ that respectively solve~\eqref{eqn:representer} and~
\eqref{eqn:second-representer-equation} are also uniquely identified. Assumption~\ref{ass:efficiency-bound}(ii) bounds $\Omega(x,z)$ away from zero and infinity. Note that conditional expectation operators can be shown to be bounded under weak conditions on the joint density \citep{carrasco2007linear}. 

\begin{corollary}\label{coro:local-efficiency}
	The   influence function \eqref{eqn:influence-function}  attains the efficiency bound of $\mu$ in $\mathcal{M}_{\np}$ at the submodel where 
	$h_0(x,z)$ solves
	\eqref{eqn:second-representer-equation}
	and
	Assumptions~\ref{ass:shadow}--\ref{ass:efficiency-bound} hold.
\end{corollary}

\subsection{Computation}

In this section, we discuss the computations of $\widehat\delta_0$ in \eqref{eqn:hat-delta0} and $\widehat h$ in~\eqref{eqn:hat-h} that are both required for $\widehat\mu_{\repa}$.
For the computation of $\widehat\delta_0$, we aim to solve the following constrained optimization problem:
\begin{equation}\label{eqn:original-optimization-for-delta0}
	\widehat\delta_0\in\argmin_{\delta\in\Delta_n} M_n(\delta),\qquad \text{s.t.} \qquad Q_n(\delta)\leq c_n.
\end{equation}
The constrained function $Q_n$ in~\eqref{eqn:sample-analogue-Q} is quadratic in $\delta$ and the objective function 
\[
M_n(\delta)=\frac{1}{n}\sum_{i=1}^n (1-R_i)\delta^2(X_i,Z_i)
\]
is also a quadratic function. 
It   remains   to impose some constraints on $\Delta_n$ so that we could have a tractable optimization problem.

A computationally simple choice for $\Delta_n$ are linear sieves as defined in~\eqref{eqn:Deltan-sieve}; that is, for any $\delta\in\Delta_n$, we have $\delta=\bm\psi^\T\bm\theta$, where $\bm\theta=(\theta_1,\ldots,\theta_{q_n})^\T$ and $\bm\psi(x,z)=\{\psi_1(x,z),\ldots,\psi_{q_n}(x,z)\}^{\T}$. For this choice of $\Delta_n$, if we define $\Delta=\{\delta:\Vert\delta\Vert_{\infty,\alpha}\leq K\}$ for some $K>0$ as required by Assumption~\ref{ass:parameter-space}, then the constraint $\Vert\bm\psi^\T\bm\theta\Vert_{\infty, \alpha}\leq K$ can be highly nonlinear in $\bm\theta$. 
We thus follow \cite{newey2003instrumental} 
by defining $\Delta$ as the closure of the set
$\{\delta:\Vert\delta\Vert_{2,\alpha_0}\leq K\}$ under $\Vert\cdot\Vert_{2,\alpha_0}$, where $\alpha_0>\alpha+(d+1)/2$ and
\begin{equation*}
	\Vert\delta\Vert_{2,\alpha_0}^2=\sum_{|\lambda|\leq \alpha_0}\int_{[0,1]^{d}\times\mathcal{Z} }\Big\{D^{\lambda}\delta(x,z)\Big\}^2dxdz,
\end{equation*}
In the above equation, $\mathcal{Z}$ denotes the support of $Z$, $\lambda$ is a $(d+1)$-dimensional vector of nonnegative integers, $|\lambda|=\sum_{i=1}^{d+1}\lambda_i$, and the operator $D^\lambda(\cdot)$ is defined in Definition~\ref{def:smooth}.
Then the constraint $\bm\psi^\T \bm\theta\in\Delta_n$ now turns  to
$\bm\theta^\T \bm H_n\bm\theta\leq K^2$,
where
\begin{equation*}
	\bm H_n=\sum_{|\lambda|\leq \alpha_0}\int_{[0,1]^{d}\times \mathcal{Z}} \Big\{D^{\lambda}\bm\psi(x,z) D^{\lambda}\bm\psi^\T(x,z)\Big\}dxdz.
\end{equation*}
Based on the above discussions, the constrained optimization problem in~\eqref{eqn:original-optimization-for-delta0} becomes:
\begin{equation*}
	\widehat{\bm\theta}_0\in\argmin_{\bm\theta} M_n(\bm\psi^\T\bm\theta),\quad \text{s.t.}\quad  Q_n(\bm\psi^\T\bm\theta)\leq c_n \quad \text{and}\quad \bm\theta^\T \bm H_n\bm\theta\leq K^2,
\end{equation*}
and the estimate of $\delta_0$ is $\widehat\delta_0=\bm\psi^\T\widehat{\bm\theta}_0$.
The computation of $\widehat h$ is similar to that of $\widehat\delta_0$. Specifically, we first solve
\begin{equation*}
	\widehat{\bm \theta}_h\in\argmin_{\bm\theta} C_n(\bm\psi^\T\bm\theta),\quad\text{s.t.}\quad \bm\theta^\T \bm H_n\bm\theta\leq K^2,
\end{equation*}
and then let $\widehat h = \bm\psi^\T\widehat{\bm \theta}_h$.
In the above optimization problems, $c_n$, $K$, and the basis function dimensions $k_n,q_n$ are all tuning parameters. Several data-driven methods for selecting tuning parameters in series estimation have been discussed in \cite{li2007nonparametric}. Here we suggest using 5-fold  cross-validation to select the tuning parameters.

	\section{Numerical studies}\label{sec:numerical-studies}

\subsection{Simulation}\label{subsec:simulation}
In this subsection, we conduct simulation studies to evaluate the performance of the proposed estimators in finite samples. We consider two different cases. In case (I), the data are generated under models where the full data distribution is identified. In  case (II), the full data distribution is not identified but Assumption~\ref{ass:representer} holds.

For  case (I), we generate  four   covariates $X=(X_1,X_2,X_3,X_4)^\T$ according to $X_j\sim U(0,1)$ for $j=1,\ldots,4$. 
We consider four data generating settings, including  combinations of  two choices of outcome models and  two choices of propensity score models. 
\begin{eqnarray*}
	&&f(Y=y\mid X=x)\thicksim\left\{\begin{array}{ll}
		N(1+2x_1+4x_2+x_3+3x_4,1),&~~~ \text{ Linear}\\
		N\{1+2x_1^2+2\exp(x_2)+\sin(x_3)+x_4, 1\}, &~~~ \text{ Nonlinear}
	\end{array}\right.\\
	&&\text{logit } f(R=1\mid x,y)=\left\{\begin{array}{ll}
		3+2x_1+x_2+x_3-0.5x_4-0.8y, & ~~~\text{ Linear}\\
		3.5+3x_1^2+4\exp(x_2)+\sin(x_3)+0.5x_4-2y,  &~~~\text{ Nonlinear}
	\end{array}\right.\\
	&&f(Z=z\mid X=x,Y=y)\sim N(3-2x_1+x_1^2+4x_2+x_3-2x_4+3y,1).
\end{eqnarray*}
The missing data proportion in each of these settings is about $50\%$.
For each setting, we replicate 1000 simulations at sample sizes $500$ and $1000$.
We apply  the proposed estimators $\widehat\mu_{\repa}$ (\REPA) and $\widehat\mu_{\rep}$ (\REP) to estimate the population outcome mean $\mu=E(Y)$. For compasion, we also
use an  inverse probability weighted  estimator (\IPW)  with a  linear-logistic propensity score model assuming MNAR  and a regression-based estimator  (\texttt{marREG}) assuming MAR to estimate $\mu$.

\begin{figure}
	\subfigure[LL]{\includegraphics[width=0.45\textwidth]{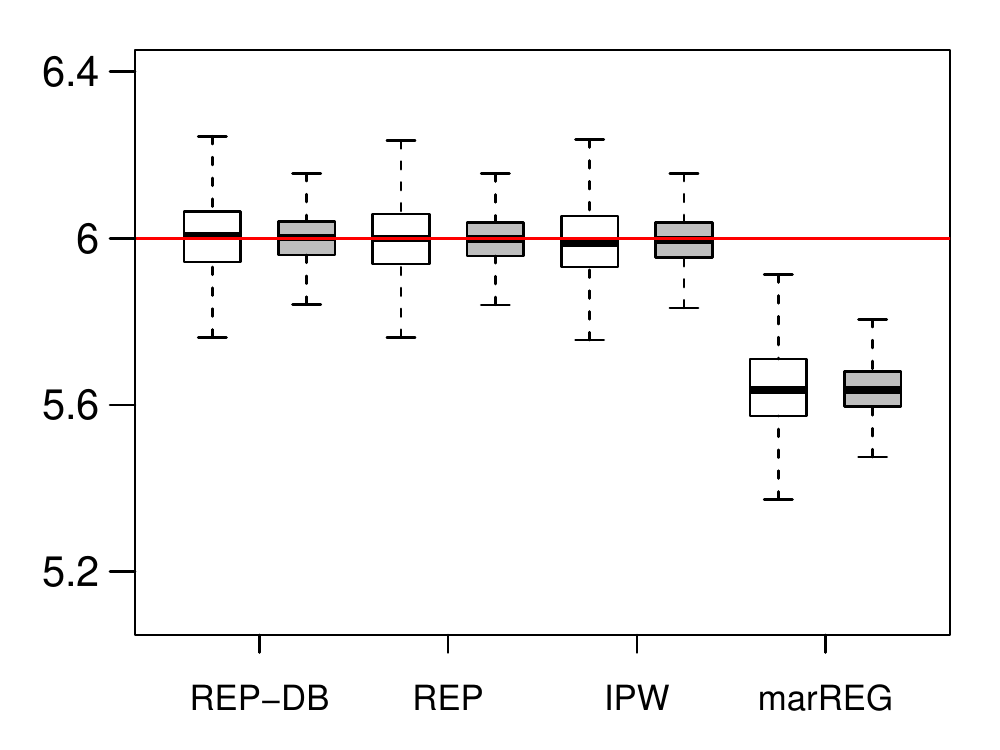}}
	\subfigure[NL]{\includegraphics[width=0.45\textwidth]{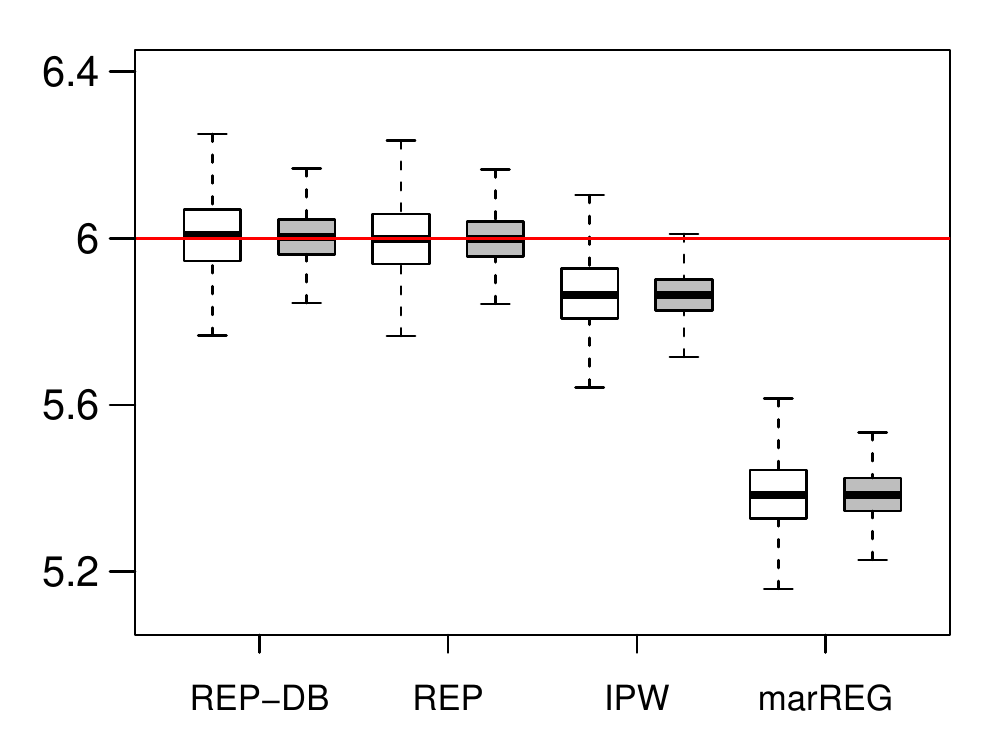}}
	\subfigure[LN]{\includegraphics[width=0.45\textwidth]{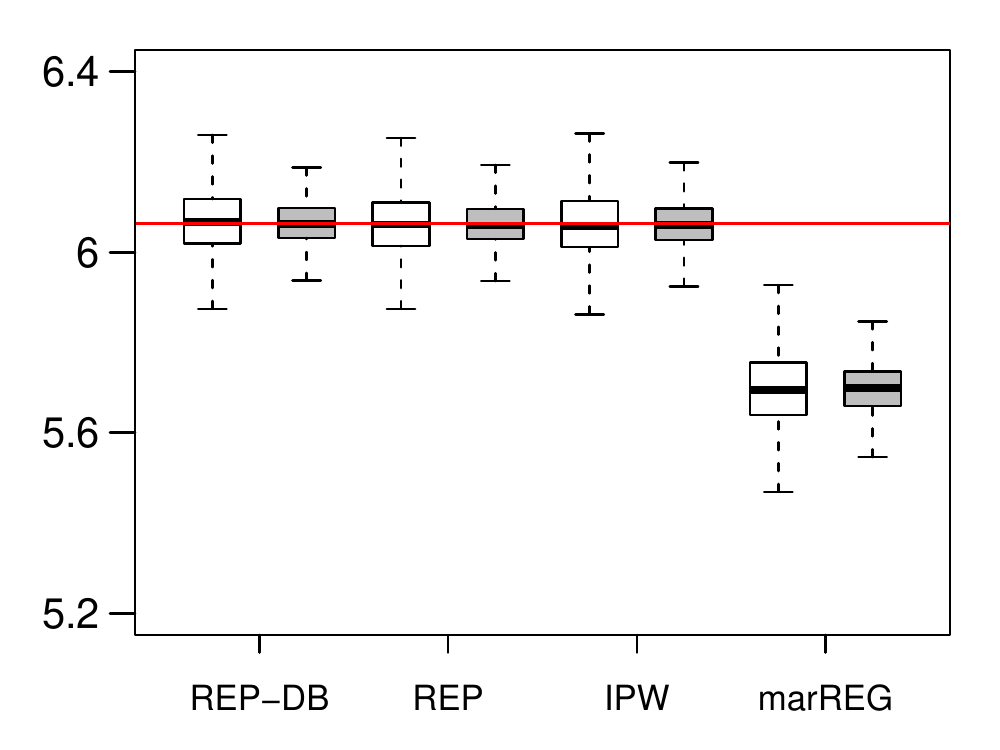}}
	\hfill
	\subfigure[NN]{\includegraphics[width=0.45\textwidth]{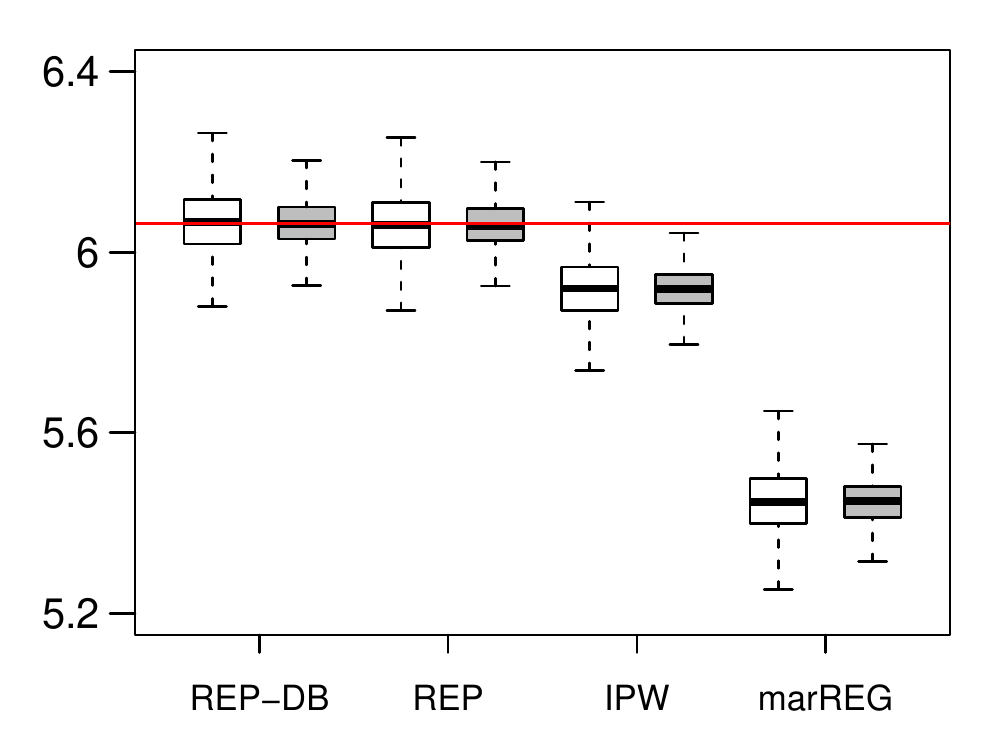}}
	\caption{Comparisons in case (I) between the proposed two estimators (\REPA~and \texttt{REP})  and existing estimators (\texttt{IPW} and \texttt{marREG}) under sample sizes $n=500$ and $n=1000$.
		The abbreviation LL stands for Linear-logistic   propensity score model  with Linear outcome model, and the other three scenarios are analogously defined.
		The horizontal line marks the true value of the outcome mean.
	}\label{fig:sim-res-boxplot}
\end{figure}

Simulation results are reported in Figure~\ref{fig:sim-res-boxplot}. 
In all four settings, the proposed estimators \REPA~and \REP~  have negligible bias. 
In contrast, the \IPW~estimator can have comparable bias with ours only when the propensity score model is correctly specified; see settings (a) and (c). 
If the propensity score model is incorrectly specified as in settings (b) and (d), the \IPW~estimator exhibits an obvious downward bias and does not vanish when the sample size increases. 
As expected, the \texttt{marREG}~estimator   has non-negligible bias in all settings.

We also calculate the 95\% confidence interval based on the proposed estimator \texttt{REP-DB} and the \IPW~estimator. 
Coverage probabilities of these two approaches  are shown in Table~\ref{tab:sim-res-cp}. 
The   \REPA~estimator based confidence intervals have coverage probabilities   close to the nominal level of $0.95$ in all scenarios even under small sample size $n=500$. 
In contrast, the  \IPW~estimator based confidence intervals have coverage probabilities  well below the nominal value
if the propensity score model is incorrectly specified.

For  case (II), we generate data according to Model 1 in Example~\ref{exam:counter}. As with case (I), we consider two different sample sizes $n=500$ and $n=1000$. We calculate the  bias (Bias), Monte Carlo standard deviation (SD) and 95\% coverage probabilities (CP) based on 1000 replications in each setting. For comparison, we also apply the \IPW~estimator with a correct propensity score model to estimate $\mu$. Since the full data distribution is not identified, the performance of \IPW~estimator depends on initial values during the optimization process. We consider two different settings of initial values for optimization parameters: true values and random values from the uniform distribution $U(0,1)$. The results are summarized in Table~\ref{tab:sim-second}.

\begin{table}
	\caption{\label{tab:sim-res-cp}Coverage probability of the 95\% confidence interval of the \REPA~and \IPW~estimators.  }
	\centering
	\begin{tabular}{lccccc}
		\addlinespace[1mm]
		\hline
		\addlinespace[1mm]
		$~n$	&Methods & LL &NL &LN &NN\\
		\hline
		\addlinespace[1.5mm]
		\multirow{ 2}{*}{~500}&\REPA &0.940 &0.932 &0.942		&0.939\\
		&\IPW &0.930 &0.635 &	0.928	&0.491 \\
		\addlinespace[1.5mm]
		\multirow{ 2}{*}{1000}&\REPA &0.945 & 0.933&	0.948	& 0.951\\
		&\IPW &0.943&0.381 &	0.951	&0.177 \\
		\addlinespace[1.5mm]
		\hline 
	\end{tabular}
\end{table}

\begin{table}
	\caption{	\label{tab:sim-second}Comparisons in case (II) between \REPA~and \texttt{IPW} under  $n=500$ and $n=1000$. }
	\centering
	\begin{tabular}{cccccccccccc}
		\hline
		\addlinespace[1mm]
		&  \multicolumn{3}{c}{\REPA} & &\multicolumn{3}{c}{\IPWt} & & \multicolumn{3}{c}{\IPWu}  \vspace{0.3mm}\\
		\cline{2-4}
		\cline{6-8}
		\cline{10-12}
		\addlinespace[1mm]
		$n$ & Bias  &SD   & CP & &Bias  &SD 
		&CP & & Bias  &SD  &CP   \\
		\addlinespace[0.8mm]
		\hline
		\addlinespace[1mm]
		\addlinespace[1mm]
		500 &0.008   &0.033   &0.917 &  &$-0.003$
		&0.050 &0.923 & & $-0.113$&0.206 &0.709  \\
		\addlinespace[1mm]
		1000 &0.003   &0.024   &0.942 &  &$-0.004$
		&0.035 &0.933 & & $-0.134$&0.216 &0.667  \\
		\addlinespace[1mm]
		\hline
	\end{tabular}
\end{table}

We observe from Table~\ref{tab:sim-second} that the proposed estimator \REPA~has negligible bias, small standard deviation and satisfactory coverage probability even under sample size $n=500$. As sample size increases to $n=1000$, the 95\% coverage probability is close to the nominal level. For the \IPW~estimator, only when the initial values for optimization parameters  are set to be true values, it has comparable performance with \REPA. However, if the initial values are randomly drawn from $U(0,1)$, the \IPW~estimator has non-negligible bias, large standard deviation and low coverage probability. As sample size increases, the situation becomes worse. We also calculate the \IPW~estimator when initial values are drawn from other distributions, e.g., standard normal distribution. The performance is even worse and we do not report the results here. The simulations in this case demonstrate the superiority of the proposed estimator over existing estimators which require identifiability of the full data distribution.

\subsection{Empirical example}

We apply the proposed methods to the China Family Panel Studies, which was previously analyzed in \cite{miao2019identification}. 
The dataset includes 3126 households in China. The outcome $Y$ is the log of current home price (in $10^4$ RMB yuan), and it has missing values due to the nonresponse of house owner and the non-availability from the real estate market.
The missingness process of home price is likely to be not at random, because subjects having expensive houses may be less likely to disclose their home prices. The missing data rate of current home price is $21.8\%$. The completely observed covariates $X$ includes 5 continuous variables: travel time to the nearest business center, house building area, family size, house story height, log of family income, and 3 discrete variables: province, urban (1 for urban househould, 0 rural), refurbish status. The shadow variable $Z$ is chosen as the construction price of a house, which is also completely observed. The construction price is related to the current price of a house, and it the shadow variable assumption that nonresponse is independent of the construction price conditional on the current price and fully observed covariates is a reasonable assumption as the construction price can be viewed as error prone proxy for the current home value, and as such is no longer predictive of the missingness mechanism once the current home value has been accounted for.

We apply the proposed  estimator \REPA~to estimate the outcome mean and the 95\% confidence interval. 
We also use the competing \IPW~estimator and two estimators assuming MAR (\texttt{marREG} and \texttt{marIPW})  for comparison. The results are shown in Table~\ref{tab:application}. We observe that the  results from the proposed estimator are similar to those from the \IPW~estimator, both yielding slightly lower estimates of home price on the log scale than those obtained from the standard MAR estimators. 
However, when the data are transformed back to the original scale, the deviations are notable and amount to approximately $1.13\times 10^4$ RMB yuan.
These analysis results are generally consistent with those in \cite{miao2019identification}. 
\begin{table}
	\caption{Point estimates and 95\% confidence intervals of the outcome mean for the home pricing example}
	\centering
	\begin{tabular}{ccccc}
		\addlinespace[1mm]
		\hline
		\addlinespace[1mm]
		Methods	& & Estimate & &95\% confidence  interval   \\
		\addlinespace[1mm]
		\hline
		\addlinespace[1mm]
		\REPA & &2.591 & &(2.520, 2.661) \\
		\addlinespace[1mm]
		\IPW & &2.611 & &(2.544, 2.678)   \\
		\addlinespace[1mm]
		\texttt{marREG} & &2.714 & &(2.661, 2.766)  \\
		\addlinespace[1mm]
		\texttt{marIPW} & &2.715 & &(2.659, 2.772)  \\
		\addlinespace[1mm]
		\hline
	\end{tabular}
	\label{tab:application}
\end{table}
	\section{Discussion}\label{sec:discussion}

With the aid of a shadow variable, we have  established a necessary and sufficient condition for nonparametric identification of mean functionals of  nonignorable missing data even if the joint distribution is not identified. Then we strengthen this condition by imposing a representer assumption that is necessary for $\sqrt{n}$-estimability of the mean functional.
The assumption
involves the existence of solutions to a representer equation, which is a Fredholm integral equation of the first kind and can be satisfied under mild requirements. Based on the representer equation, we propose a  sieve-based estimator for the mean functional, which bypasses the difficulties of correctly specifying and estimating the unknown missingness mechanism and the outcome regression. Although the joint distribution is not identifiable, the proposed estimator is shown to be consistent for the  mean functional. 
In addition, we establish conditions under which the proposed estimator is    asymptotically normal.  We  would like to point out that since the solution to the representer equation is not uniquely determined, one cannot simply apply standard theories for nonparametric sieve estimators to derive the above asymptotic results. In fact, we need to first construct a consistent estimator for the solution set, and then find from the estimated set a consistent estimator for an appropriately chosen solution.	We finally show that the proposed estimator
attains the semiparametric efficiency bound for the shadow variable model at a key submodel where the representer is uniquely identified. 

The availability of a valid shadow variable is crucial for the proposed approach. 
Although it is generally not possible to test the shadow variable assumption via observed data without making another untestable assumption, the existence of such a variable is practically reasonable in  the empirical example presented in this paper and similar situations where one or more proxies or surrogates of a variable prone to missing data may be available. In fact, it is not uncommon in survey studies and/or cohort studies in the health and social sciences, that certain outcomes may be sensitive and/or expensive to measure accurately, so that a gold standard measurement is obtained only for a select subset of the sample, while one or more proxies or surrogate measures may be available for the remaining sample. Instead of a standard measurement error model often used in such settings which requires stringent identifying conditions, the more flexible shadow variable approach proposed in this paper provides a more robust alternative to incorporate surrogate measurement in a nonparametric framework, under minimal identification conditions. Still, the validity of the shadow variable assumptions generally requires domain-specific knowledge of experts and needs to be investigated on a case-by-case basis.
As advocated by \cite{robins2000sensitivity}, in principle, one can also conduct sensitivity analysis to assess how results would change if the shadow variable assumption were violated by some pre-specified amount. 

The proposed methods may be improved or extended in several directions. Firstly, the proposed identification and estimation framework may be extended to handle nonignorable missing outcome regression or missing covariate problems. Secondly, one can use modern machine learning techniques to solve the representer equation so that an improved estimator  may be achieved that adapts to sparsity structures in the data. Thirdly, it is  of great interest to extend our results to handling other  problems of   coarsened data, for instance,   unmeasured confounding problems in causal inference.  
We plan to pursue these and other related issues  in future research.

\section*{Appendix}\label{appn} 


\subsection*{Existence of solutions to representer equation~\eqref{eqn:representer} under completeness conditions}

We adopt the singular value decomposition (\cite{carrasco2007linear}, Theorem 2.41) of compact operators to characterize conditions for existence of a solution to~\eqref{eqn:representer}. Let $L^2\{F(t)\}$ denote the space of all square-integrable functions of $t$ with respect to a cumulative distribution function $F(t)$, which is a Hilbert space with inner product $\langle g,h\rangle=\int g(t)h(t)dF(t)$. Let $K_x$ denote the conditional expectation operator $L^2\{F(z\mid x,r=1)\}\rightarrow L^2\{F(y\mid x,r=1)\}$, $K_xh=E\{h(Z)\mid x,y,r=1\}$ for $h\in L^2\{F(z\mid x,r=1)\}$, and let $(\lambda_n,\varphi_n,\psi_n)_{n=1}^{+\infty}$ denote a singular value decomposition of $K_x$. We assume the following regularity conditions:
\begin{manualcond}{A1}\label{cond:A1}
	$\int \int f(z\mid x, y,r=1)f(y\mid x, z,r=1)dydz<+\infty$
\end{manualcond}
\begin{manualcond}{A2}\label{cond:A2}
	$\int \tau^2(x,y)f(y\mid x,r=1)dy<+\infty$.
\end{manualcond}
\begin{manualcond}{A3}\label{cond:A3}
	$\sum_{n=1}^{+\infty}\lambda_{n}^{-2}\big|\langle \tau(x,y),\psi_n\rangle \big|^2<+\infty$.
\end{manualcond}
Given $f(z\mid x,y,r=1)$, the solution to~\eqref{eqn:representer} must exist if the completeness condition and Conditions \ref{cond:A1}--\ref{cond:A3} all hold. The proof follows immediately from Picard's theorem (\cite{kress1989linear}, Theorem 15.18) and Lemma 2 of \cite{miao2018identify}.

\section*{Supporting information}
\label{SM}

Supporting information includes additional lemmas and proofs of all the theoretical results.

		\bibliographystyle{apalike}
	\bibliography{mybib}
\end{document}